\documentclass[10pt]{article}
\usepackage{amsmath,amsthm,amsfonts,amssymb}
\usepackage{graphicx}
\usepackage{hyperref}
\usepackage[usenames]{color}

\newtheorem{theorem}{Theorem}[section]
\newtheorem{lemma}[theorem]{Lemma}
\theoremstyle{definition}

\newtheorem{corollary}[theorem]{Corollary}

\newtheorem{definition}[theorem]{Definition}
\newtheorem{example}[theorem]{Example}

\newtheorem{remark}[theorem]{Remark}
\newtheorem{proposition}[theorem]{Proposition}
\newcounter{comcount}

\renewcommand{\L}{{\mathcal{L}}}

\def\A{{\mathcal{A}}}
\def\B{{\mathcal{B}}}
\def\Ce{{\mathcal{C}}}
\def\M{{\mathcal{M}}}
\def\N{{\mathcal{N}}}

\def\E{{\mathcal{E}}}
\def\T{{\mathcal{T}}}

\def\CC{{\mathcal C}}
\def\CF{{\mathcal F}}
\def\CP{{\mathcal P}}

\newcommand{\Th}{{\mathrm{Th}}}
\newcommand{\Tr}{{\mathrm{T}}}
\newcommand{\At}{{\mathrm{At}}}
\newcommand{\Mod}{{\mathrm{Mod}}}
\newcommand{\Diag}{{\mathrm{Diag}}}
\newcommand{\Con}{{\mathrm{Con}}}
\newcommand{\atp}{{\mathrm{atp}}}

\newcommand{\K}{{\mathbf{K}}}
\renewcommand{\H}{{\mathbf{H}}}
\renewcommand{\S}{{\mathbf{S}}}
\renewcommand{\P}{{\mathbf{P}}}
\newcommand{\Pw}{{\mathbf{P_{\!\! \omega}}}}
\newcommand{\Ps}{{\mathbf{P_{\!\! s}}}}
\newcommand{\Pf}{{\mathbf{P_{\! f}}}}
\newcommand{\Pu}{{\mathbf{P_{\!\! u}}}}

\newcommand{\Ld}{{\underrightarrow{\mathbf{L}}}}
\newcommand{\Ls}{{\underrightarrow{\mathbf{L}}_{\! s}}}
\newcommand{\Loc}{{\mathbf{L}}}
\newcommand{\e}{{\mathbf{\,_e}}}
\newcommand{\Op}{{\mathbf{O}}}

\newcommand{\pvar}{{\mathbf{Pvar}}}
\newcommand{\var}{{\mathbf{Var}}}
\newcommand{\qvar}{{\mathbf{Qvar}}}
\newcommand{\ucl}{{\mathbf{Ucl}}}
\newcommand{\Res}{{\mathbf{Res}}}
\newcommand{\Dis}{{\mathbf{Dis}}}
\newcommand{\Cat}{{\mathbf{Cat}}}

\newcommand{\Rad}{{\mathrm{Rad}}}
\newcommand{\V}{{\mathrm{V}}}

\newcommand{\AT}{{\mathfrak{A}}}
\newcommand{\BT}{{\mathfrak{B}}}
\newcommand{\TT}{{\mathfrak{T}}}

\title{Unification theorems in algebraic geometry}

\author{E.\,Daniyarova, A.\,Myasnikov, V.\,Remeslennikov}

\date{August 12, 2008}

\begin{document}
\maketitle

\begin{abstract}
In this paper, for a given finitely generated algebra (an
algebraic structure with arbitrary operations and no predicates)
$\A$ we study finitely generated limit algebras of $\A$, approaching
them via model theory and algebraic geometry. Along the way we lay down 
foundations of algebraic geometry over arbitrary algebraic structures.
\end{abstract}

\tableofcontents

\section{Introduction}
\label{se:intro}

Quite often relations between sets of elements of a fixed
algebraic structure $\A$ can be described in terms  of equations
over $\A$. In the classical case, when $\A$ is a field, the area
of mathematics where such relations are studied is known under the
name of {\it algebraic geometry}. It is natural to use the same name
in the general case. Algebraic geometry over arbitrary algebraic
structures is a new  area of research in modern algebra,
nevertheless, there are already several breakthrough particular results here,
as well as, interesting developments of a general theory.
 Research in this area started with a series
of papers by Plotkin \cite{Plot1,Plot2}, Baumslag, Kharlampovich,
Myasnikov, and Remeslennikov \cite{BMR1,MR2, KM1,KM2}.

There are general results which hold in the algebraic geometries
over arbitrary algebraic structures, we refer to them as the  {\it
universal algebraic geometry}.  The main purpose of this paper is
to lay down the basics of the universal algebraic geometry in a
coherent form. We emphasize here the relations between model
theory, universal algebra, and algebraic geometry. Another goal is
quite  pragmatic~---  we intend to unify here some common methods
known in different fields under different names. Also, there are
several essentially the same results that independently occur in
various branches of modern algebra, were they  are treated by
means specific to the area. Here we give very general proofs of
these results based on  model theory and universal algebra.

Limit algebras, in all their various incarnations,  are the  main
object of this paper. The original notion came from group theory
where limit groups play a prominent part. The limit groups of a
fixed group $G$ appear in many different situations: in
combinatorial group theory as groups discriminated by $G$
($\omega$-residually $G$-groups or fully residually $G$-groups)
\cite{bb,gb,MR2,BMR2,BMR3}, in the algebraic geometry over groups
as the coordinate groups of irreducible varieties over $G$
\cite{BMR1,KM1,KM2,KM3,Sela1}, groups universally equivalent to
$G$ \cite{Remesl,GS,MR2}, limit groups of $G$ in the
Grigorchuk-Gromov's metric \cite{CG}, in the theory of equations
in groups \cite{Lyndon, Razborov1, Razborov2,
KM1,KM2,KM3,Groves2}, in group actions
\cite{BF,GLP,MRS,Groves1,G1}, in the solutions of Tarski problems
\cite{KM4,Sela2}, etc. These numerous characterizations of limit
groups make them into a very robust tool linking  group theory,
topology and logic. It turned out that many of the results on
limit groups can be naturally generalized to Lie algebras
\cite{Daniyarova1,DKR1,DKR2,Daniyarova2,Daniyarova3}.

Our prime objective is to convey some basic facts of the general
theory of limit  algebras in an arbitrary language. We prove the
so-called unification theorems for limit groups that show that  the
characterization results above hold in the general case as well.


\section{Preliminaries}
\label{se:preliminary}

\subsection{Languages and structures}
\label{sec:formulas}

  Let $\L = \CF \cup \CP \cup \CC$ be a first-order  {\em language}
  (or a {\em signature}), consisting of a set $\CF$ of symbols
  of operations $F$ (given together with their arities $n_F$), a
  set $\CP$ of symbols of predicates $P$ (given together with their arities
  $n_P$) and a set of constants $\CC$. If $\CP = \emptyset$ then the language
  $\L$  is {\em functional},  whereas  $\L$ is {\em relational} if $\CF =\CC = \emptyset$.

  For languages $\L_1 \subseteq \L_2$ we say that $\L_1$ is a {\em reduct} of $\L_2$
  and  $\L_2$ is an {\em expansion} of $\L_1$. The language $\L^{fun} = \L \smallsetminus{\CP}$
  is the functional part of $\L$. From now we fix a first-order functional language $\L$.
  Almost everything we prove holds (under appropriate adjustments) for arbitrary languages,
  but the exposition for functional languages is shorter.

\begin{example}
The language of {\em groups} consists of a binary operation
$\cdot$ (multiplication), a unary operation $^{-1}$ (inversion),
and a constant symbol $e$ or $1$ (the identity).
\end{example}

\begin{example}
The language of {\em unitary rings} consists of three binary
operations $+$, $-$ and $\cdot$ (addition, subtraction and
multiplication), and constants $0$ and $1$.
\end{example}

An $\L$-{\em structure} $\M$ is given by the following data: (i) a
non-empty set $M$ called the {\em universe} of $\M$; (ii) a
function $F^{\M}:M^{n_F}\to M$ of arity $n_F$ for each $F \in
\CF$; (iii) an element $c^{\M}\in M$ for each $c\in \CC$. We often
white the structure as $\M = \langle \, M; \, F^{\M}, \, c^{\M},
\, F\in \CF, \, c\in \CC\,\rangle$. We refer to $F^\M$ and $c^\M$
as {\em interpretations} of the symbols $F$ and $c$ in $\M$, and
sometimes omit superscripts $\M$ (when the interpretation is
obvious from the context). Typically we denote structures in $\L$
by capital calligraphic letters and their universes (the
underlying sets) by the corresponding capital Latin letters.
Structures in a functional language are termed {\em algebras} (or
{\em universal algebras}). An algebra $\E$ with the universe
consisting of a single element is called {\em trivial}. Obviously,
interpretation of symbols from $\L$ in $\E$ is unique.

As usual, one can define the notion of a  homomorphism, and all
its variations, between structures in a given language. If a
subset $N \subseteq M$ is closed under the operations $F$ of $\CF$
and contains all the constants $c \in \CC$ then restrictions of
the operations $F$ onto $N$, together with the constants $c$,
determine a new $\L$-structure, called a {\em substructure} $\N$
of $\M$, in which case we write $\N \leq \M$. For  a subset
$M^\prime  \subseteq M$ the intersection of all substructures of
$\M$ containing $M^\prime $ is a substructure $\M^\prime $ of $\M$
generated by $M^\prime $ (so $M^\prime $ is a generating set for
$\M^\prime$), symbolically $\M^\prime  = \langle M^\prime
\rangle$. $\M$ is termed {\em finitely generated} if it has a
finite generating set.

Let $X = \{x_1, x_2, \ldots \}$ be a finite or countable set of
variables.
  {\em Terms}
  in $\L$ in variables $X$ are formal expressions defined recursively as follows:
   \begin{itemize}
   \item  [T1)] variables $x_1, x_2, \ldots, x_n, \ldots $ are terms;
   \item [T2)] constants from $\L$ are terms;
   \item [T3)] \label{def:term} if $F(x_1, \ldots, x_n) \in \CF$
   and $t_1, \ldots, t_n$ are terms then $F(t_1, \ldots, t_n)$ is a term.
  \end{itemize}
For $F \in \CF$ we write $F(x_1, \ldots,x_n)$ to indicate that $n
= n_F$.

By $\Tr_\L = \Tr_\L(X)$ we denote the set of all terms in $\L$.
For a term $t \in \Tr_\L$ one can define the set of variables
$V(t) \subset X$ that occur in $t$.  We write $t(x_1, \ldots,x_n)$
to indicate that $V(t) \subseteq \{x_1, \ldots,x_n\}$. Also, we
use the vector notation  $t(\bar{x})$, where
$\bar{x}=(x_1,\ldots,x_n)$. Following the recursive definition of
$t$ one can define  in a natural way a  function $t^{\M}: M^n \to
M$ (which we sometimes again denote by $t$). If $V(t) = \emptyset$
then $t$ is a {\em closed} term and $t^\M$ is just a constant.
Observe, that the universe of the substructure of $\M$ generated
by a subset $M^\prime \subseteq M$ is equal to $\bigcup \{t(M')
\mid t\in \Tr_{\L}(X)\}$, where $t(M^\prime)$ is the range of the
function $t$.

The condition T3) allows one to define an operation $F^{\T_{\L}
(X)}$ on the set of terms $\Tr _{\L}(X)$. By T2) the set $\Tr
_{\L}(X)$ contains all constants from $\L$, which gives a natural
interpretation of constants in $\Tr _{\L}(X)$. These altogether
turn the set $\Tr _{\L}(X)$ into  an $\L$-structure $\T _{\L}(X)$,
which is called the {\em absolutely free} $\L$-algebra with basis
$X$. The name comes from the the following universal property of
$\T _{\L}(X)$: for any $\L$-structure $\M$ a map $h: X \to M$,
extends to a unique $\L$-homomorphism $h: \T _{\L}(X) \to \M$.

{\em Formulas} in $\L$ (in variables $X$) are defined recursively
as follows:
 \begin{itemize}
   \item  [F1)] if $t, s \in T_\L(X)$ then $(t = s)$ is a formula (called an {\em atomic formula});
   \item [F2)] if $\phi$ and $\psi$ are formulas then
   $\neg \phi, (\phi \vee \psi), (\phi \wedge \psi), (\phi \rightarrow \psi)$  are formulas;
   \item [F3)] If $\phi$ is a formula and $x$ is a variable then $\forall x \phi$ and $\exists x \phi$ are formulas.
  \end{itemize}

 For a formula $\phi$ one can define the set $V(\phi)$
 of {\em free variables} of $\phi$ according to the rules F1)--F3).
 Namely, $V(t_1 = t_2) = V(t_1) \cup V(t_2)$, $V(\neg \phi) = V(\phi)$,
 $V(\phi \circ \psi) = V(\phi) \cup V(\psi)$, where  $\circ \in \{\vee, \wedge, \rightarrow\}$,
 and $V(\forall x \phi) = V(\exists x \phi) = V(\phi) \smallsetminus \{x\}$.
 We write $\phi(x_1, \ldots, x_n)$ in the case when $V(\phi) \subseteq \{x_1, \ldots, x_n\}$.
 Let $\Phi _{\L} (X)$ be the set of all formulas in $\L$ with $V(\phi) \subseteq X$.
 A formula $\phi$ with $V(\phi) = \emptyset$  termed  a {\em sentence}, or a {\em closed} formula.

 If $\phi(x_1, \ldots, x_n)  \in \Phi _{\L} (X)$
 and $m_1, \ldots,m_n \in M$ then one can define,  following the conditions F1)--F3),
 the relation ``$\phi$ is true  in $\M$  under the interpretation  $x_1 \to m_1, \ldots, x_n \to
 m_n$''
 (symbolically $\M \models \phi(m_1, \ldots,m_n)$).
 It is convenient sometimes to view this relation as an
 $n$-ary predicate $\phi^M$  on $M$. If $h:X \to M$ is
 an interpretation of variables then we denote  $\phi^h = \phi^M(h(x_1), \ldots, h(x_n))$.

 A set of formulas $\Phi \subseteq \Phi _{\L} (X)$
 is {\em consistent} if there is an $\L$-structure $\M$ and
 an interpretation $h:X \to M$ such that $\M \models \phi^h$
 for   every $\phi \in \Phi$.   In this case one says that $\Phi$ is {\em realized} in $\M$.

The following result is due to Malcev, it plays a crucial role in
model theory.

\bigskip
\noindent {\bf Theorem} [Compactness Theorem] {\em  Let $\K$ be a
class of $\L$-structures and $\Phi \subseteq \Phi _{\L} (X)$. If
every finite subset of $\Phi$ is realized in some structure in
$\K$ then the whole set $\Phi$ is realized in some ultraproduct of
structures from $\K$.}

 \bigskip


\subsection{Theories}

 Two formulas $\phi, \psi \in \Phi _{\L} (X)$ are called {\em equivalent}
 if $\phi^h = \psi^h$ for any interpretation $h:X \to M$ and any $\L$-structure $\M$.
 One of the principle results in mathematical logic states that any formula $\phi \in \Phi _{\L} (X)$
 is equivalent to a formula $\psi$ in the following form:
 \begin{equation}
 \label{eq:prenex}
 Q_1 x_1 \ldots Q_m x_m \: \left ( \bigvee_{i = 1}^n (\bigwedge_{j= 1}^k \psi_{ij}) \right ),
  \end{equation}
 where $Q_i \in \{\forall, \exists\}$ and $\psi_{ij}$  is an atomic formula  or its negation.
 One of the standard ways to characterize  complexity of formulas is according to their  quantifier prefix
 $ Q_1 x_1 \ldots Q_m x_m$ in (\ref{eq:prenex}).

 If in (\ref{eq:prenex}) all the quantifiers $Q_i$ are universal
 then the formula $\psi$ is called {\em universal} or $\forall$-formula,
 and if all of them are existential then $\psi$ is   {\em existential}
 or $\exists$-formula.  In this fashion $\psi$ is $\forall\exists$-formula
 if the prefix has only one alteration of quantifiers (from $\forall$ to $\exists$).
 Similarly, one can define $\exists\forall$-formulas.
 Observe, that $\forall$- and $\exists$-formulas are dual relative to negation,
 i.e.,  the negation of $\forall$-formula is equivalent to an $\exists$-formula,
 and the negation of $\exists$-formula is equivalent to an $\forall$-formula.
 A similar result holds for $\forall\exists$- and $\exists\forall$-formulas.
 One may consider formulas with more alterations of quantifiers,
 but we have no use of them in this paper.

 A formula in the form (\ref{eq:prenex}) is {\em positive}
 if it does not contain negations (i.e., all $\psi_{ij}$ are atomic).
 A formula is {\em quantifier-free} if it does not contain quantifiers.
 We denote the set of all quantifier-free formulas from  $\Phi _{\L} (X)$
 by $\Phi_{\mathrm{qf},\L}(X)$, and the set of all atomic formulas   by  $\At_\L(X)$.

 Recall that a {\em theory} in the language $\L$ is
 an arbitrary consistent set of sentences in $\L$. A theory $T$ is {\em complete}
 if for every sentence $\phi$ either $\phi$ or $\neg\phi$ lies in $T$. By $\Mod (T)$
 we denote the (non-empty) class of all $\L$-structures  $\M$ which satisfy all
 the sentences from $T$. Structures from $\Mod(T)$ are termed {\em models} of $T$ and
 $T$ is a set of {\em axioms} for the class $\Mod(T)$. Conversely, if $\K$
 is a class of $\L$-structures then the set $\Th(\K)$ of sentences, which are
 true in all structures from $\K$,  is called the {\em elementary theory} of $\K$.
 Similarly, the set $\Th_\forall(\K)$ ($\Th_\exists(\K)$) of all $\forall$-sentences
 ($\exists$-sentences) from $\Th(\K)$ is called the {\em universal} ({\em existential})
 theory of $\K$. The following notions play an important part  in this paper.
 Two $\L$-structures $\M$ and $\N$  are  \textit{elementarily equivalent}
 if  $\Th (\M)=\Th (\N)$, and they are {\em universally} ({\em existentially})
 equivalent if $\Th _{\forall}(\M)=\Th _{\forall}(\N)$  ($\Th
 _{\exists} (\M)=\Th _{\exists}(\N)$). In this event we write,
 correspondingly, $\M\equiv \N$, $\M\equiv_{\forall} \N$
 or $\M\equiv_{\exists} \N$. Notice, that due to the duality
 mentioned above $\M\equiv_{\forall} \N \Longleftrightarrow \M\equiv_{\exists} \N$
 for arbitrary $\L$-structures $\M$ and $\N$.

A class of $\L$-structures $\K$ is {\em axiomatizable} if $\K =
\Mod(T)$ for some theory $T$ in $\L$. In particular, $\K$ is
$\forall$- ($\exists$-, or $\forall\exists$-) axiomatizable if the
theory $T$ is $\forall$- ($\exists$-, or $\forall\exists$-)
theory.


\section{Algebras}
\label{subsec:algebras}

There are several types of classes of $\L$-structures that play a
part in general algebraic geometry: prevariaeties, quasivarieties,
universal closures, and $\A$-algebras. We refer to \cite{MR2} for
a detailed discussion on this and related matters. Here we present
only a few  properties and characterizations of these classes,
that will be used in the sequel.  Most of them are known and can
be found in the classical books on universal algebra, for example,
in \cite{Malcev}. On the algebraic theory of quasivarieties, the
main subject of this section, we refer to \cite{Gorbunov}.


\subsection{Congruences}

In this section we remind some notions and introduce notation on
presentation of algebras via generators and relations.

Let $\M$ be an arbitrary fixed $\L$-structure. An equivalence
relation $\theta$ on  $M$  is a {\em congruence} on $\M$  if for
every operation $F \in \CF$ and any elements $m_1,\ldots,m_{n_{\!
F}}$, $m'_1,\ldots,m'_{n_{\! F}}\in M$ such that $m_i\sim
_{\theta} m'_i$, $i=1, \ldots, n_{\! F}$, one has
$F^{\M}(m_1,\ldots,m_{n_{\! F}})\sim _{\theta}
F^{\M}(m'_1,\ldots,m'_{n_{\! F}})$.

For a congruence $\theta$ the operations $F^\M$, $F \in \CF$,
naturally induce well-defined operations on the factor-set
$M/\theta$. Namely, if  we denote by $m/\theta$ the equivalence
class of $m \in M$ then $F^{\M/\theta}$ is defined by
$$
F^{\M/\theta}(m_1/\theta,\ldots,m_{n_{\!
F}}/\theta)=F^{\M}(m_1,\ldots,m_{n_{\! F}})/\theta
$$
 for any $m_1,\ldots,m_{n_{\! F}} \in M$. Similarly, $c^{\M/\theta}$
 is defined  for $c \in \CC$ as the class $c^\M/\theta$.
 This turns the factor-set $\M/\theta$
 into an $\L$-structure.
 It follows immediately from the construction that the map
$h: M \to M/\theta$, such that $h(m)=m/\theta$, is an
$\L$-epimorphism $h: \M \to \M/\theta$, called the {\em canonical}
epimorphism.

The set $\Con (\M)$ of all congruences on $\M$ forms a lattice
relative to the inclusion $\theta_1\leqslant \theta_2$, i.e.,
every two congruences in $\Con (\M)$ have the least upper and the
greatest lower  bounds in the ordered set $\langle \Con(\M),
\leqslant \rangle$. To see this, observe first that the
intersection of an arbitrary set $\Theta = \{\theta_i, i \in I\}$
of congruences on $\M$ is again a congruence on $\M$, hence the
greatest lower bound for $\Theta$. Now, the intersection of the
non-empty set $\{\theta \in \Con(\M) \mid \theta_i \leqslant
\theta \ \forall \: \theta_i \in \Theta\}$ is the least upper
bound for $\Theta$. The following result is easy.

\begin{lemma}\label{rem}
Let $\M$ be an $\L$-algebra, $\{\theta_i \mid i \in I\}\subseteq
\Con (\M)$  and $\theta=\bigcap_{i \in I} \theta_i$. Then
$\M/\theta$ embeds into the direct product  $\prod_{i\in
I}{\M/\theta_i}$ via the diagonal monomorphism $m/\theta \to
\prod_{i\in I}{m/\theta_i}$.
\end{lemma}

A homomorphism $h:\M \to \N$ of two $\L$-structures determines the
{\em kernel} congruence $\ker h$ on $\M$, which is defined by
$$
m_1\sim_{\ker h} m_2 \quad \Longleftrightarrow \quad
h(m_1)=h(m_2), \quad m_1,m_2 \in M.
$$
Observe, that if  $\theta\in \Con (\M)$ and  $\theta \leqslant
\ker h$  then the map  $\bar{h}:\M/\theta \to \N$ defined by
$\bar{h} (m/\theta)=h(m)$ for  $m\in M$ is a homomorphism of
$\L$-structures.

\begin{definition}
\label{de:congruent} A set of atomic formulas $\Delta \subseteq
\At _{\L}(X)$ is called {\em congruent} if the binary relation
$\theta_\Delta$ on the set of terms $\Tr _{\L}(X)$ defined by
(where $t_1, t_2 \in \Tr _{\L}(X)$)
$$
t_1\sim_{\theta_{\Delta}} t_2 \quad \Longleftrightarrow \quad
(t_1=t_2)\in \Delta.
$$
is a congruence on the free $\L$-algebra $\T _{\L}(X)$.
\end{definition}

The following lemma characterizes congruent sets of formulas.

\begin{lemma}
\label{le:congruent} A set of atomic formulas $\Delta \subseteq
\At _{\L}(X)$ is congruent if and only if it satisfies the
following conditions:
\begin{enumerate}
\item $(t=t)\in \Delta$ for any term $t\in \Tr _{\L}(X)$;
\item if $(t_1=t_2) \in \Delta$ then  $(t_2=t_1) \in \Delta$ for any terms $t_1, t_2 \in \Tr _{\L}(X)$;
\item if $(t_1=t_2) \in \Delta$ and  $(t_2=t_3) \in \Delta$
then $(t_1=t_3) \in \Delta$ for any terms  $t_1, t_2, t_3 \in
\Tr_{\L} (X)$;
\item if $(t_1=s_1),\ldots, (t_{n_{F}}=s_{n_{F}}) \in \Delta$ then
$(F(t_1,\ldots,t_{n_{F}})=F(s_1,\ldots,s_{n_{F}})) \in \Delta$ for
any terms  $t_i,s_i \in \Tr_{\L} (X)$, $i=1, \ldots,n_{F}$, and
any functional symbol $F \in \L$.
\end{enumerate}
\end{lemma}

  \begin{proof}
  Straightforward.
  \end{proof}

Since the intersection of an arbitrary set of congruent sets of
atomic formulas is again congruent, it follows that for a set
$\Delta \subseteq \At _{\L}(X)$ there is the least congruent
subset $[\Delta] \subseteq \At _{\L}(X)$, containing $\Delta$.
Therefore, $\Delta$ uniquely determines the congruence
${\theta_{\Delta}} = {\theta_{[\Delta]}}$.

For an $\L$-algebra  $\M$ generated by a set $M^\prime \subseteq
M$ put  $X=\{x_m \mid  m\in M'\}$ and consider a set
$\Delta_{M^\prime}$ of all atomic formulas $(t_1=t_2) \in \At_\L
(X)$ such that $\M\models (t_1=t_2)$ under the interpretation $x_m
\to m, m\in M^\prime$.  Obviously, $\Delta_{M^\prime}$ is a
congruent set in $\At_\L (X)$ (the set of all relation in $\M$
relative to $M^\prime$). A subset $S \subseteq \Delta_{M^\prime}$
is called a set of {\em defining relations} of $\M$ relative to
$M^\prime$ if $[S] = \Delta_{M^\prime}$.  In this event the pair
$\langle X \mid S\rangle$ termed a {\em presentation} of $\M$ by
generators $X$ and relations $S$.

\begin{lemma}\label{xs}
If  $\langle X \mid S \rangle$ is a presentation of $\M$ then  $\M
\cong \T_{\L}(X)/\theta_{S}$.
\end{lemma}

\begin{proof}
The map  $h':X\to M'$ defined by $h'(x_m)= m$, $m\in M'$, extends
to a homomorphism $h:\T_{\L}(X) \to \M$. Clearly, $t_1\sim_{\ker
h} t_2$ if and only if $(t_1=t_2) \in [S]$  for terms  $t_1,t_2
\in \Tr_{\L}(X)$. Therefore, $\T_{\L}(X)/\theta_{S} \cong
\T_{\L}(X)/\ker h$. Now the result follows from the isomorphism
$\T_{\L}(X)/\ker h \cong \M$.
\end{proof}


\subsection{Quasivarieties}
\label{subsec:quasivar}

In this section we discuss quasivarieties and related objects. The
main focus is on how to generate the least quasivariety containing
a given class of structures $\K$. A model example here is the
celebrated Birkhoff's theorem which describes $\var(\K)$, the
smallest variety containing $\K$, as the class $\H \S \P (\K)$
obtained from $\K$ by taking direct products (the operator $\P$),
then substructures (the operator $\S$), and then homomorphic
images (the operator $\H$). Along the way we introduce some other
relevant operators. On the algebraic theory of quasivarieties we
refer to \cite{Gorbunov} and \cite{Malcev}.

We fix, as before, a functional language $\L$ and a class of
$\L$-algebras $\K$. We always assume that $\K$ is an abstract
class, i.e., with  any algebra $\M \in \K$ the class $\K$ contains
all isomorphic copies of $\M$.

Recall that an {\em identity} in $\L$ is a formula of the type
$$
\forall x_1 \ldots \forall  x_n \left
(t(x_1,\ldots,x_n)=s(x_1,\ldots,x_n)\right ),
$$
where $t,s$ are terms in $\L$. Meanwhile, a  {\em quasi-identity}
is a formula of the type
$$
\forall x_1 \ldots \forall x_n \left( (\bigwedge\limits_{i=1}^m
t_i (\bar {x})=s_i (\bar {x})) \quad \rightarrow \quad (t(\bar
{x})=s(\bar {x})) \right),
$$
where $t(\bar{x}), s(\bar{x}), t_i(\bar{x}), s_i(\bar{x})$ are
terms in $\L$ in variables $\bar{x} = (x_1,\ldots,x_n)$.

A class of $\L$-structures is  called a {\em quasivariety} ({\em
variety}) if it can be axiomatized by a set of quasi-identities
(identities).  Given a class of $\L$-structures $\K$ one can
define the quasivariety $\qvar(\K)$, generated by $\K$, as the
quasivariety axiomatized by the set $\Th_{\rm qi}(\K)$ of all
quasi-identities which are true in all structures from $\K$, i.e.,
$\qvar(K) = \Mod(\Th_{\rm qi}(\K))$. Notice, that $\qvar(\K)$ is
the least quasivariety containing $\K$. Similarly, one defines the
variety $\var(\K)$  generated by $\K$.

Observe, that an identity  $\forall\, \bar {x} (t(\bar {x})=s(\bar
{x}))$ is equivalent to a quasi-identity $\forall\,   \bar {x}
(x=x \rightarrow t(\bar {x})=s(\bar {x}))$, therefore, $\qvar(\K)
\subseteq \var(\K)$.

Before we proceed with quasivarieties, we introduce one more class
of structures. Namely, $\K$ termed  a {\em prevariety} if $\K =
\S\P(\K)$.  By $\pvar(\K)$ we denote the least prevariety,
containing $\K$.  The prevariety $\pvar(\K)$  grasps the residual
properties of the structures from $\K$. An $\L$-structure $\M$ is
{\em separated} by $\K$ if for any pair of non-equal elements
$m_1,m_2\in M$ there is a structure $\N \in \K$ and a homomorphism
$h: \M \to \N$ such that $h(m_1)\neq h(m_2)$. By $\Res (\K)$ we
denote the class of $\L$-structures separated by $\K$.

In the following lemma we collect some known facts on
prevarieties.

\begin{lemma}
\label{le:prevar} For any class of $\L$-structures $\K$ the
following holds:
    \begin{enumerate}
    \item [1)] \label{item:1} $\pvar (\K)=\S \P (\K) \subseteq \qvar (\K)$;
    \item [2)] \label{item:2} $\pvar (\K)=\Res (\K)$;
    \item [3)] \label{item:3} $\pvar(\K)$ is axiomatizable if and only if $\pvar(\K) = \qvar(\K)$.
    \end{enumerate}
\end{lemma}

\begin{proof}
Equality 1) follows directly from definitions.

2) was proven for groups in \cite{MR2}, here we give a general
argument. It is easy to see that $\Res (\K)$ is a prevariety, so
$\pvar (\K) \subseteq \Res (\K)$. To show converse, take a
structure  $\M \in \Res (\K)$ and consider the  set  $I$ of all
pairs $(m_1,m_2)$, $m_1,m_2\in M$, such that $m_1\neq m_2$. Then
for every  $i\in I$ there exists a structure $\N_i \in \K$ and a
homomorphism $h_i:\M \to \N_i$ with $h_i(m_1)\ne h_i(m_2)$. The
homomorphisms $h_i, i\in I,$ give rise to the ``diagonal''
homomorphism $h : \M \to \prod_{i\in I} {\N_i}$, which is
injective by construction. Hence $\M \in \S \P (\K)$, as required.

3) is due to Malcev \cite{Malcev2}.
\end{proof}

Prevarieties play an important role in combinatorial algebra, they
can be characterized as classes of structures admitting
presentations by generators and relator. Namely, let $X$ be a set
and $\Delta$ a set of atomic formulas from $\Phi_\L(X)$. Following
Malcev \cite{Malcev}, we say that a presentation $\langle X \mid
\Delta\rangle$ {\em defines  a structure} $\M$ in a class $\K$ if
there is a map $h:X \to M$ such that
 \begin{enumerate}
 \item [D1)] $h(X)$ generates $\M$ and all the formulas from $\Delta$
 are realized in $\M$ under the interpretation $h$;
 \item [D2)] for any structure $\N \in \K$ and any map $f:X \to N$
 if all the formulas from $\Delta$ are realized in $\N$ under $f$
 then there exists a unique homomorphism $g:\M \to \N$ such that $g(h(x)) = f(x)$ for every $x \in X$.
 \end{enumerate}
If $\langle X \mid \Delta\rangle$ defines a structure in $\K$ then
this structure is unique up to isomorphism, we denote it by
$F_\K(X,\Delta)$.

\bigskip
{\bf Theorem} \cite{Malcev} {\it A class $\K$, containing the
trivial system $\E$,  is a prevariety if and only if any
presentation $\langle X \mid \Delta\rangle$ defines a structure in
$\K$.}

\bigskip
To present  similar characterizations for quasivarieties we need
to introduce the following operators.

As was mentioned above,  $\P (\K)$ is the class of direct products
of structures from $\K$. Recall, that the  direct product of
$\L$-structures  $\M_i$, $i\in I$, is an $\L$-structure  $\M =
\prod_{i\in I}{\M_i}$ with the universe $M=\prod_{i\in I}{M_i}$
where the functions and constants from $\L$ are interpreted
coordinate-wise. If all the structures $\M_i$ are isomorphic to
some structure $\N$ then  we refer to $\prod_{i\in I}{\M_i}$ as to
a direct power of $\N$ and denote it by $\N^{I}$. By  $\Pw (\K)$
we denote the class of all finite direct products of structures
from $\K$.

Recall, that a substructure $\N$ of a direct product $\prod_{i\in
I}{\M_i}$ is a {\em subdirect product} of the structures  $\M_i$,
$i\in I$, if $p_j(\N)=\M_j$ for the canonical projections $p_j:
\prod_{i\in I}{\M_i} \to M_j$,  $j\in I$. By  $\Ps (\K)$ we denote
the class of all subdirect products  of structures from $\K$.

Let  $I$ be a set, $D$  a filter over $I$ (i.e., a collection $D$
of subsets of $I$ closed under finite intersections and such that
if $a \in D$ then $b \in D$ for any $b \subseteq I$ with $a
\subseteq b$, and also we assume that $\emptyset \not \in D$), and
$\{M_i \mid i\in I\}$ a family of sets. On the direct product
$\prod_{i\in I}{M_i}$ one can define an equivalence relation $\sim
_D$ such that $a \sim _D b$  if and only if $\{i\in I \mid p_i(a)
= p_i(b)\} \in D$. We denote the factor-set by $\prod_{i\in
I}{M_i}/D$, and the equivalence class of an element $a$ by $a/D$.
Now, if $\{\M_i \mid i\in I\}$  is  a collection of
$\L$-structures then the equivalence $\sim _D$ becomes a
congruence on the direct product $\prod_{i\in I}{\M_i}$, in which
case  the filterproduct $\M = \prod_{i\in I}{\M_i}/D$ of the
structures $\M_i, i \in I,$ over $D$ is defined as the
factor-structure $\prod_{i\in I}{\M_i}/\sim_D$. If $D$ is an
ultrafilter on $I$ (a filter that contain either $a$ or $I
\smallsetminus a$ for any $a \subseteq I$) then a  filterproduct
over $D$ is called an {\em ultraproduct}, furthermore, if all the
structures $\M_i$ are isomorphic to some structure $\N$ then the
ultraproduct $\prod_{i\in I}{\M_i}/D$ is called an {\em
ultrapower} and we denote it by $\N^I/D$. By  $\Pf (\K)$ and $\Pu
(\K)$  we denote, correspondingly, the classes of filterproducts
and ultraproducts of structures from $\K$.

 Let  $\K_e = \K \cup \{\E\}$, where $\E$ is the trivial $\L$-structure
 introduced earlier.  A word of warning is needed here. Sometimes, direct
 products $\prod_{i\in I}{\M_i}$ are defined being equal to $\E$  for the
 empty set $I$ (see, for example, \cite{Gorbunov}), but we elect not to do so,
 assuming always that $I$ is non-empty  and adding $\E$ to the class,  if needed.

\begin{lemma}
For any class of $\L$-structures $\K$ the following holds:
\begin{enumerate}
\item  [5)] \label{st4} $\qvar (\K)=\S \Pf (\K)\e$;
\item [6)] \label{st5} $\qvar (\K)=\S \P \Pu (\K)\e = \S \Pu \P (\K)\e$;
\item [7)] \label{st6} $\qvar (\K)=\S \Pu \Pw (\K)\e$;
\end{enumerate}
\end{lemma}

\begin{proof}
 5) is due to Malcev  \cite[\textsection 11, Theorem~4]{Malcev}.
6) and 7) are due to Gorbunov \cite[Corollary~2.3.4,
Theorem~2.3.6]{Gorbunov}.
\end{proof}

Now we give another characterization of quasivarieties, for this
we need to introduce direct limits.

Recall, that a partial ordering  $(I, \leqslant)$ is {\em
directed} if any two elements from $I$ have an upper bound.  A
triple  $\Lambda=(I, \M_i, h_{ij})$, consisting of a directed
ordering $(I, \leqslant)$, a set of $\L$-structures $\{\M_i, i\in
I\}$, and a set of homomorphisms $h_{ij}:\M_i \to \M_j$ ($i,j \in
I$, $i \leqslant j$), is called a {\em direct system} of
structures $\M_i$, $i\in I$, if
\begin{enumerate}
\item $h_{ii}$ is the identity map for every $i\in I$;
\item $h_{jk}\circ h_{ij}=h_{ik}$ for any $i,j,k \in I$ with $i\leqslant j\leqslant k$.
\end{enumerate}
We call a directed system $\Lambda=(I, \M_i, h_{ij})$ {\em
epimorphic} if all the homomorphisms $h_{ij}:\M_i \to \M_j$ are
surjective.

Given a direct system $\Lambda=(I, \M_i, h_{ij})$ one can consider
an equivalence relation $\equiv$ on a set $\{(m_i,i) \mid  m_i\in
M_i, \: i\in I\}$ defined by
$$
(m_i,i)\equiv (m_j,j) \quad \Leftrightarrow \quad \exists \: k\in
I, \; i,j\leqslant k,\quad h_{ik} (m_i)=h_{jk}(m_j).
$$
By $\langle m, i\rangle$ we denote the equivalence class of
$(m,i)$ under $\equiv$. Now one can turn the factor-set
 $M = \{(m_i,i) \mid  m_i\in M_i, \: i\in I\}/\equiv$ into an $\L$-structure
 $\M$ interpreting the constants and functions from $\L$ as follows:
\begin{enumerate}
\item if $c\in \L$ is a constant then
$c^{\M}=\langle c^{\M_i},i \rangle$ for an arbitrary chosen $i\in
I$;
\item if $F\in \L$ is a function and
$\langle m_1,i_1 \rangle,\ldots, \langle m_{n_{\! F}},i_{n_{\! F}}
\rangle \in M$ then
$$F^{\M}(\langle m_1,i_1 \rangle,\ldots, \langle m_{n_{\!
F}},i_{n_{\! F}} \rangle)= F^{\M_j}(\langle h_{i_1 j}(m_1),i_1
\rangle,\ldots, \langle h_{i_{n_{\! F}} j}(m_{n_{\! F}}),i_{n_{\!
F}} \rangle)$$ for an arbitrary chosen $j\in I$ with  $i_1,\ldots,
i_{n_{\! F}}\leqslant j$.
\end{enumerate}
The structure $\M$ is well-defined, it is called the {\em direct
limit} of the system $\Lambda$, we denote it by
$\underrightarrow{\lim} \M_i$. It is easy to see that
$\underrightarrow{\lim} \M_i$ has the following property. Let
$i\in I$ be a fixed index. Put $J_i=\{j\in I \mid  i \leqslant
j\}$ a nd denote $\Lambda_i=(J_i, \M_j, h_{jk}, \: j,k\in J_i)$.
Then $\Lambda_i$ is a direct system whose direct limit $\M^i$ is
isomorphic to $\M$. By $\Ld (\K)$ and $\Ls (\K)$ we denote the
class of direct and epimorphic direct limits of structures from
$\K$.

The following result gives a characterization of quasivarieties in
terms of direct limits.

\begin{lemma}
For any class of $\L$-structures $\K$ the following holds:
$$\qvar (\K)=\S \Ls \P (\K)\e = \Ls \S \P (\K)\e = \Ls \Ps (\K)\e = \Ld \S \P (\K)\e.$$
\end{lemma}
 \begin{proof}
 See \cite[Corollary~2.3.4]{Gorbunov}.
 \end{proof}


\subsection{Universal closures}
\label{subsec:ucl}

In this section we study the universal closure $\ucl(\K) = \Mod
(\Th _{\forall} (\K))$ of a given class of $\L$-structures $\K$.

Structures from $\ucl(\K)$ are determined by local properties of
structures from $\K$. To explain precisely we need to introduce
two more operators.

 Recall \cite{BMR2,MR2}, that a structure $\M$ is {\em discriminated}
 by $\K$ if  for any finite set $W$ of elements from $\M$ there is a structure $\N \in \K$ and a
 homomorphism $h :\M \to \N$ whose restriction onto $W$ is injective. Let $\Dis (\K)$ be the class
 of $\L$-structures  discriminated by $\K$. Clearly, $\Dis(\K) \subseteq \Res(\K)$.

   To introduce the second operator we  need to describe  local submodels of a structure $\M$.
   First, we replace the language $\L$ by a new relational language $\L^{rel}$, where
   every operational and constant symbols $F \in \CF$ and $c \in \CC$ are replaced, correspondingly,
   by a new predicate symbol $R_F$ of arity $n_F +1$  and  a new unary predicate symbol $R_c$.
   Secondly, the structure $\M$ turns into a $\L^{rel}$-structure $\M^{rel}$, where
   the predicates $R_c^{\M^{rel}}$ and $R_F^{\M^{rel}}$ are defined by
 \begin{enumerate}
  \item [R1)] for $m \in M$ the predicate  $R_c^{\M^{rel}}(m)$ is true in $\M^{rel}$ if and only if $c^{\M}=m$;
  \item [R2)] for $m_0,m_1,\ldots,m_{n_{F}}\in M$ the predicate
  $R_F^{\M^{rel}}(m_0,m_1,\ldots,m_{n_{F}})$ is true in $\M^{rel}$
  if and only if $F^{\M}(m_1,\ldots,m_{n_{F}})=m_0$.
 \end{enumerate}

 Third, if $\L_0$ is a finite  reduct  (sublanguage) of $\L$
 then by $\M^{\L_0}$ we denote the reduct of $\M^{rel}$,
 where only predicates corresponding to constants and
 operations  from  $\L_0$ are survived,  so $\M^{\L_0} $ is an
 $\L^{rel}_0$-structure.
 Now, following \cite{Malcev}, by a {\em local submodel} of $\M$ we
 understand a finite substructure of $\M^{\L_0}$ for some finite reduct  $\L_0$ of $\L$.

 Finally, a structure $\M$ is {\em locally embeddable} into $\K$
 if every local submodel of $\M$ is isomorphic to some local submodel of
 a structure from $\K$ (in the language $\L^{rel}_0$). By $\Loc (\K)$ we
 denote the class of $\L$-structures locally embeddable into $\K$.

It is convenient for us to rephrase the notion of a  local
submodel in terms of formulas.

Let  $\L^\prime$ be a finite reduct of $\L$ and $X$ a finite set
of variables.  A  quantifier-free formula  $\varphi$ in
$\L^\prime$  is called a {\em diagram-formula} if $\varphi$ is a
conjunction  of atomic formulas or their negations that satisfies
the following conditions:
 \begin{enumerate}
\item [1)] every formula $\neg(x = y)$,  for each pair  $(x, y) \in X^2$ with
$x \neq y$, occurs in $\varphi$;
\item [2)] for each functional symbol  $F \in \L^\prime$ and each
tuple of variables $(x_0,x_1,\ldots,x_{n_{F}}) \in X^{n_F+1}$
either formula $F(x_1,\ldots,x_{n_{F}})=x_0$ or its negation
occurs in $\varphi$;
\item [3)] for each constant symbol  $c\in\L'$ and each  $x \in X$
either $x= c$ or its negation $\neg(x=c)$  occurs in $\varphi$.
\end{enumerate}

We say that $\varphi$ is a {\em diagram-formula} in $\L$ if it is
a diagram-formula for some finite reduct $\L^\prime$ of $\L$ and a
finite set $X$. The name of diagram-formulas comes from the
diagrams of algebraic structures (see Section
\ref{subsec:algebras}).

The following lemma is easy.
\begin{lemma}
 For any local submodel $\N$ of $\M$ there is a diagram-formula $\varphi_\N(X)$
 in a finite set of variables $X$ of cardinality $|N|$ such that $\M \models \varphi_\N(h(X))$
 for some bijection $h :X \to N$. And conversely, if  $\M \models \varphi(h(X))$ for
 some diagram-formula $\varphi(X)$ in $\L$ and an interpretation $h:X \to M$ then
 there is a local submodel $\N$ of $\M$ with the universe $h(X)$ such that $\varphi = \varphi_\N$
 (up to a permutation of conjuncts).
\end{lemma}

\begin{corollary} \label{co:diagram-form}
An $\L$-structure $\M$ is locally embeddable into a class $\K$ if
and only if every diagram-formula realizable in $\M$ is realizable
also in some structure from $\K$.
\end{corollary}

\begin{lemma} \label{le:ucl}
For any class of $\L$-structures $\K$ the following holds:
\item [8)] $\ucl (\K)=\Loc (\K)$;
\item [9)] $\ucl (\K)=\S \Pu (\K)$;
\item [10)] $\Dis (\K) \subseteq \ucl (\K)$;
\item [11)] $\Ld (\K) \subseteq \ucl (\K)$.
\end{lemma}
   \begin{proof}
   To prove 8) and 9) we show that $\Loc (\K)\subseteq \S \Pu (\K) \subseteq \ucl (\K)
\subseteq \Loc (\K)$. The first inclusion has been proven by
Malcev \cite{Malcev}, but we briefly discuss it for the sake of
completeness. Let $\M$ be a structure from $\Loc (\K)$. By
Corollary \ref{co:diagram-form} every diagram-formula $\varphi$
realizable in $\M$ is realizable also in some structure
$\N_\varphi$ from $\K$. By the Compactness Theorem
 the set $\Phi_\M$ of all diagram-formulas realizable in $\M$ is
 realized in some ultraproduct $\N = \prod_{\varphi}\N_\varphi/D$,
 where $\varphi$ runs over $\Phi_\M$. By Lemma \ref{le:diag-core}
 the core $\Diag_0(\M)$  of the diagram of $\M$  is also realized
 in $\N$ under an appropriate interpretation of constants
 $c_m, m \in M$ (see Section  \ref{subsec:algebras}).
 Now the substructure of $\N$ generated by all elements
 $c_m, m \in M$,  is isomorphic to $\M$. Hence $\M \in \S \Pu (\K)$.

Inclusion $\S \Pu (\K) \subseteq \ucl (\K)$ follows from two known
results: any universal class is closed under substructures (which
is obvious) and the Los theorem ~\cite{Malcev,Marker}. To see that
$\ucl (\K) \subseteq \Loc (\K)$ consider an arbitrary $\M \in \ucl
(\M)$. If $\varphi (x_1,\ldots,x_n)$ is a diagram-formula which is
realized in $\M$ then a universal sentence $\psi=\forall
x_1,\ldots, x_n \neg \, \varphi (x_1,\ldots,x_n)$ is false in
$\M$. Hence, there exists a structure $\N \in \K$ on which  $\psi$
is false, so $\N \models \neg \psi$. Therefore, $\varphi
(x_1,\ldots,x_n)$ is realized in $\N$. By Corollary
\ref{co:diagram-form} $\M \in \Loc (\K)$, as required.

To see 10) it suffices to notice that $\Dis (\K) \subseteq \Loc
(\K)$ and then apply 8).

11) follows from 9) and \cite{Gorbunov} (Theorem 1.2.9), where it
is shown that  $\Ld (\K) \subseteq \S \Pu (\K)$.

   \end{proof}


\subsection{$\A$-Algebras} \nopagebreak
\label{subsec:algebras} \label{sec3}

Let $\A$ be a fixed $\L$-algebra. In this section we discuss
$\A$-algebras~--- principal objects in algebraic geometry over
$\A$. Informally, an $\A$-algebra is an $\L$-algebra with a
distinguished subalgebra $\A$. Even though this notion seems
simple, one needs to develop a formal framework to deal with
$\A$-algebras. It will be convenient to use two equivalent
approaches: one is categorical and another is  logical (or
axiomatic).

\begin{definition} \label{de:cat-A} [Categorical]
An $\A$-algebra is a pair $(\B,\lambda)$, where $\B$ is an
$\L$-algebra and $\lambda :\A\to \B$  is an embedding.
  \end{definition}

For the  axiomatic definition we are going to use the language of
diagrams. By $\L_\A$ we denote the language $\L \cup \{c_a \mid a
\in A\}$, which is obtained from $\L$ by adding a new constant
$c_a$ for every element $a \in A$.

Observe, that every $\A$-algebra $(\B,\lambda)$ can be viewed as
an $\L_\A$-algebra when the constant $c_a$ is interpreted by
$\lambda(a)$.

Recall that by $\At_{\L_\A}(\emptyset)$ we denote the set of all
atomic sentences in the language $\L_\A$. The {\em diagram} $\Diag
(\A)$ of $\A$ is the set of all  atomic sentences from $
\At_{\L_\A}(\emptyset)$  or their negations which are true in
$\A$. To work with diagrams we need to define  several related
sets of formulas.

The {\em core} $\Diag_0(\A)$ of the diagram $\Diag(\A)$ consists
of the following formulas:
\begin{itemize}
\item $c = c_a$ for each constant symbol $c \in \CC$ and $a \in A$ such that $c^\A = a$;
\item $F(c_{a_1},\ldots, c_{a_{n_{F}}})=c_{a_0}$,
for each functional symbol $F \in \CF$ and each tuple of elements
$(a_0, a_1,\ldots, a_{n_{F}}) \in A^{n_F+1}$ such that
$F^{\A}(a_1,\ldots, a_{n_{F}})=a_0$;
\item $c_{a_1}\neq c_{a_2}$, for each pair $(a_1, a_2) \in A^2$ such that $a_1 \neq a_2$.
\end{itemize}

The following result is easy
\begin{lemma}
  \label{le:diag-core}
  For an $\L$-algebra $\A$ the following hold:
   \begin{enumerate}
    \item [C1)] For every $\L_\A$-structure $\B$ if $\B \models \Diag_0(\A)$ then $\B \models \Diag(\A)$;
    \item [C2)] If $S$ is a finite subset of $\Diag_0(\A)$
    then there is a diagram-formula $\varphi(X)$ in $\L$
    and an interpretation $h: X \to A$ such that every formula from $S$ occurs
    as a conjunct in $\varphi(h(X))$ (after replacing $h(x)$ with $c_{h(x)}$)  and $\A \models
    \varphi(h(X))$;
    \item [C3)] If $\varphi(X)$ is a diagram-formula in $\L$ and $h:X \to A$
    is an interpretation such that $\A \models \varphi(h(X))$ then every
    conjunct of $\varphi(X)$ (where $x$ is replaced with $c_{h(x)}$)  belongs to $\Diag(\A)$.
    \end{enumerate}
\end{lemma}

The following result gives an axiomatic way to describe
$\A$-algebras.

\begin{lemma} \label{le:axiomatic}
 Let $\B$ be an $\L$-algebra and $\lambda: A \rightarrow B$ a map.
 Then  $(\B,\lambda)$ is an $\A$-algebra if and only if $\B \models \Diag(\A)$,
 where  $c_a$ is interpreted by $\lambda(a)$ for every $a \in A$.
 \end{lemma}
   \begin{proof} Straightforward.

   \end{proof}

This leads to the following, equivalent, definition of
$\A$-algebras.

\begin{definition}\label{de:cat-A} [Axiomatic]
An  algebra $\B$ in the language $\L_\A$ is called an $\A$-algebra
if $\B \models \Diag(\A)$.
 \end{definition}
Put
$$
\Diag ^{+} (\A) = \{\varphi \in\At_{\L_\A}(\emptyset) \mid
\A\models \varphi \},
$$
$$
\Diag ^{-} (\A) = \Diag(\A) \smallsetminus \Diag ^{+} (\A).
$$

Let $\Cat (\A)$ be the class of all $\A$-algebras. Since
$\A$-algebras are $\L_\A$-structures the standard notions of a
$\L_\A$-homomorphism, $\L_\A$-substructure, $\L_\A$-generating
set, etc.,  are defined in $\Cat(\A)$. Sometimes, we refer to them
as to an  $\A$-homomorphism, $\A$-substructure, $\A$-generating
set, etc. Class $\Cat (\A)$ with $\A$-homomorphism forms a  category of $\A$-algebras.

All the operators $\Op$ introduced in Sections \ref{subsec:ucl}
and \ref{subsec:quasivar} are defined for $\L_\A$-structures, but,
a priori, the resulting $\L_\A$-algebra may not be in the class
$\Cat (\A)$. Nevertheless, one can check directly for each such
operator $\Op$ (with the exception of the operator $\K \to \K_e$
that adds the trivial structure $\E$ to $\K$)   that
$\Op(\Cat(\A)) \subseteq \Cat(\A)$.  Sometimes, we add the
subscript $\A$ and write $\Op_\A$ to emphasize the fact that the
algebras under consideration are $\A$-algebras. Another, shorter,
way to prove this is  to show that $\Cat(\A)_e$ is a quasivariety,
and then these results, as well as some others, will follow for
free.

\begin{lemma}\label{le:A-quasivar}
The class  $\Cat (\A)\e$ is a quasivariety in the language
$\L_{\A}$ defined by the following set of quasi-identities:
\begin{enumerate}
\item $c=c_a$, for each constant symbol $c\in\L$
and element $a \in A$ such that $a=c^{\A}$;
\item $F (c_{a_1},\ldots, c_{a_{n_{F}}})=c_a$,
for each functional symbol $F \in \L$ and each tuple $(a_1,\ldots,
a_{n_{F}},a)\in A^{n_F +1}$ such that $F^{\A}(a_1,\ldots,
a_{n_{F}}) = a$;
\item $\forall \,x\: \forall\, y\: (c_{a_1}=c_{a_2} \to x=y)$,
for each pair of elements $a_1, a_2\in A$ with $a_1 \neq a_2$.
\end{enumerate}

\end{lemma}

\begin{proof}
It is easy to see that any $\A$-algebra and the trivial algebra
$\E$ satisfy the formulas above.  One needs to check the converse.
Suppose $\Ce$ is an  $\L_{\A}$-algebra, satisfying the formulas
above. If $\Ce = \E$ then $\Ce\in\Cat (\A)\e$. Assume now that
$\Ce \neq \E$. The formulas 1) and 2) show that  $\Ce \models
\Diag_0(\A) \cap \Diag^{+}(\A)$, while the formulas 3) provide
$\Ce \models \Diag_0(\A) \cap \Diag^{-}(\A)$. Altogether, $\Ce
\models \Diag_0(\A)$, so by Lemma \ref{le:diag-core}  $\Ce \models
\Diag(\A)$, as claimed.
\end{proof}

\begin{corollary}\label{oper2} Let $\A$ be an algebra and $\K$ a class of $\A$-algebras. Then the following holds:
\begin{enumerate}
\item [1)]  $\K$
is closed under the operators $\S_{\A}$,
$\P_{\A}$, $\Pw_{\A}$, $\Ps_{\A}$, $\Pf_{\A}$, $\Pu_{\A}$,
$\Ld_{\A}$, $\Ls_{\A}$, $\Loc_{\A}$;
 \item [2)] every algebra in the classes $\pvar _{\A} (\K)$,  $\ucl _{\A}
(\K)$,  $\Res _{\A}(\K)$, and  $\Dis _{\A}(\K)$ is an
$\A$-algebras;

\item [3)] every algebra  in $\qvar _{\A} (\K)$, with the exception of $\E$, is an  $\A$-algebra.
    \end{enumerate}
\end{corollary}


\section {Types, Zariski topology, and coordinate algebras}\label{sec:limit} \nopagebreak

In this section we introduce algebras defined by complete atomic
types.

\subsection{Quantifier-free types and Zariski topology}

Let $\L$ be a functional language, $T$ a theory in $\L$, and $X =
\{x_1, \ldots,x_n\}$ a finite set of variables. Recall (see, for
example, \cite{Marker}), that a {\em type} in variables $X$ of
$\L$ over $T$ is a consistent with $T$ set $p$ of formulas in
$\Phi_\L(X)$, i.e, a subset $p \subseteq \Phi_\L(X)$ that can be
realized in a structure from $\Mod(T)$.

A type $p$ is {\em complete} if it is a maximal type in
$\Phi_\L(X)$ with respect to inclusion. It is easy to see that if
$p$ is a maximal type in $X$ then for every formula $\varphi \in
\Phi _{\L}(X)$ either $\varphi \in p$ or $\neg\,\varphi \in p$.

\begin{definition}
A set $p$ of atomic or negations of atomic formulas from
$\Phi_\L(X)$  is called an {\em atomic type} in $X$ relative to a
theory $T$ if $p \cup T$ is consistent. A maximal atomic type in
$\Phi_\L(X)$  with respect to inclusion termed a {\em complete
atomic type} of $T$.
\end{definition}
  It is not hard to see that if $p$ is a complete atomic type
  then for every atomic formula  $\varphi \in \At_{\L}(X)$
  either $\varphi \in p$ or $\neg\,\varphi \in p$.

\begin{example}
Let $\M$ be an $\L$-structure and $\bar{m}=(m_1,\ldots,m_n) \in
M^n$. Then the set $\atp^{\M}(\bar{m})$ of atomic or negations of
atomic formulas from $\Phi_\L(X)$ that are true in $\M$ under an
interpretation
 $x_i \mapsto m_i$, $i=1, \ldots,n$, is a complete atomic type
 relative to any theory $T$ such that $\M \in \Mod(T)$.
\end{example}

We say that a complete atomic type $p$ in variables $X$ is
realized in $\M$ if $p = \atp^{\M}(\bar{m})$  for some $ \bar{m}
\in M^n$.

Every type $p$ in $T$ can be realized in some model of $T$ (i.e.,
a structure from $\Mod(T)$). If $p$ cannot be realized in a
structure $\M$ then we say that $\M$ {\em omits} $p$. There are
deep results in model theory on how to construct models of $T$
omitting a given type or a set of types.

For an atomic type $ p \subseteq  \Phi_\L(X)$ by $p^{+}$ and
$p^{-}$ we denote, correspondingly,  the set of all atomic and
negations of atomic formulas in $p$.

If $S$ is a set of atomic formulas from $\Phi_\L(X)$ and $\M$ is
an $\L$-structure then by $\V_\M(S)$ we denote the set $\{(m_1,
\ldots,m_n) \in M^n \mid \M \models S(m_1, \ldots,m_n)\}$ of all
tuples in $M^n$ that satisfy all the formulas from $S$. The set
$\V_\M(S)$ is called the {\em algebraic set} defined by $S$ in
$\M$.
We refer to $S$  as  a {\em system of equations} in $\L$, and to elements
of $S$ - as  {\em equations} in $\L$. Sometimes, to emphasize  that formulas
are from $\L$ we call such equations (and systems of equations) {\em coefficient-free equations},
meanwhile,  in the case  when $\L = \L_\A$, we refer to such equations as
equations with coefficients in algebra $\A$.

Following \cite{BMR1} we define  {\em Zariski topology} on $M^n, n
\geq 1,$ where algebraic sets form a prebasis of closed sets,
i.e., closed sets in this topology are obtained from the algebraic
sets by finite unions and (arbitrary) intersections.

If $p$ is an atomic type in $\L$ in variables $X = \{x_1,
\ldots,x_n\}$ then  $\V_\M(p^+)$ is an algebraic set in $M^n$.
More generally, for an arbitrary type $p$ in $X$ by $p^+$ we
denote the set of all positive formulas in $p$, i.e., all formulas
in the prenex form that do not have the negation symbol.

  If $p$ is quantifier-free type, i.e., a type consisting of
  quantifier-free formulas, then formulas in  $p^+$ are conjunctions
  and disjunctions of atomic formulas.

\begin{lemma}
Let $\M$ be an $\L$-structure and $n \in \mathbb{N}$. Then for a
subset $V \subseteq M^n$ the following conditions are equivalent:
 \begin{itemize}
 \item $V$ is closed in the Zariski topology on $M^n$;
 \item $V = \V_\M(p^+)$ for some quantifier-free type $p$ in variables  $\{x_1,
 \ldots,x_n\}$.
 \end{itemize}
 Here $\V_\M(p^+)=\{(m_1,
\ldots,m_n) \in M^n \mid \M \models p^{+}(m_1, \ldots,m_n)\}$.
\end{lemma}
\begin{proof}
Straightforward.
\end{proof}


\subsection{Coordinate algebras and complete types}

Let $\M$ be an $\L$-algebra. For a set $S$ of atomic formulas from
$\Phi_\L(X)$ denote by $\Rad_\M(S)$ the set of all atomic formulas
from $\Phi_\L(X)$ that hold on every tuple from $\V_\M(S)$. In
particular, if $\V_\M(S) = \emptyset$ then  $\Rad_\M(S) =
\At_\L(X)$. It is not hard to see that $\Rad_\M(S)$ is a congruent
set of formulas, hence it defines a congruence that we denote by
$\theta_{\Rad (S)}$. The $\L$-structure $\T_\L(X)/\theta_{\Rad
(S)}$ is called the {\em coordinate algebra} of the algebraic  set
$\V_\M(S)$. If $Y = \V_\M(S)$ then the coordinate algebra
$\T_\L(X)/\theta_{\Rad (S)}$ is denoted by $\Gamma (Y)$ and $\Rad
(S)$~ -- by $\Rad (Y)$.

The following result gives a characterization of the coordinate
algebras over an algebra $\M$.

\begin{proposition}\label{res}
A finitely generated $\L$-algebra  $\Ce$ is the coordinate algebra
of some non-empty algebraic set over an $\L$-algebra $\M$ if and
only if $\Ce$ is separated by $\M$.
\end{proposition}

\begin{proof}
Let $Y$ be an algebraic set in $M^n$. With a point
$p=(m_1,\ldots,m_n) \in M^n$ we associate a homomorphism $h_p:
\T_{\L}(X)\to \M$ defined by  $h_p(t)=t^{\M}(m_1,\ldots,m_n)$.
Clearly,
$$
\theta_{\Rad (Y)}= \bigcap\limits_{p\in Y} {\ker h_p}.
$$
Therefore, the diagonal homomorphism $\prod_{p \in
Y}:\T_{\L}(X)\to \prod_{p\in Y} \M$ induces a monomorphism
$$\Gamma (Y) = \T_{\L}(X)/\theta_{\Rad(Y)} \to \M^{|Y|}.$$
It follows that $\Gamma (Y) \in  \S\P(\M)$. Now, by Lemma
\ref{le:prevar}  $\S\P(\M)=\Res(\M)$, so $\Gamma (Y) \in \Res
(\M)$.

Suppose now that  $\Ce$ is a finitely generated $\L$-algebra from
$\Res (\M)$ with a finite generating set $X= \{x_1, \ldots,
x_n\}$.  Let $\Ce = \langle X \mid  S\rangle$ be a presentation of
$\Ce$ by the generators $X$ and relations  $S \subseteq
\At_{\L}(X)$. In this case  $\Ce$ is isomorphic to
$\T_{\L}(X)/\theta_S$. To prove that $\Ce$ is the coordinate
algebra of some algebraic set over $\M$ it suffices to show that
$\Rad _{\M}(S)=[S]$. If $(t_1=t_2)\not\in [S]$ then there exists a
homomorphism $h:\Ce \to \M$ with
$t^{\M}_1(h(x_1),\ldots,h(x_n))\ne
t^{\M}_2(h(x_1),\ldots,h(x_n))$. Obviously,
$(h(x_1),\ldots,h(x_n)) \in \V _{\M}(S)$ so $(t_1=t_2)\not\in \Rad
_{\M}(S)$. This shows that $\Rad _{\M}(S)=[S]$.
\end{proof}

\begin{lemma}
 Let $p$ be a complete atomic type in variables $X$. Then:
  \begin{itemize}
  \item   $p^{+}$  is a congruent set of formulas;
  \item $p^+ = \Rad_\M(p^+)$ for every $\L$-structure  $\M$ with $\V_\M(p) \neq \emptyset$.
  \end{itemize}
 \end{lemma}

 \begin{proof}
 Indeed, since $p$ is realized in some model $\M$ of $T$ its
 positive part $p^+$  satisfies the assumptions of Lemma \ref{le:congruent},
 hence it is congruent. It follows that $p^+$ determines a
 congruence $\theta_p$ on $\T _{\L}(X)$. Since $p$ is complete one has $p^+ = \Rad_\M(p^+)$.
\end{proof}

\begin{definition}
Let $X$ be a finite set of variables and $p$ a complete atomic
type in variables $X$.  Then the factor-algebra $\T
_{\L}(X)/\theta_p$ of the free $\L$-algebra $\T _{\L}(X)$ is
termed the algebra defined by the type $p$ and the tuple
 $(x_1/\theta_p,\ldots,x_n/\theta_p)$  is called a generic point of $p$.
\end{definition}

Clearly, any complete atomic type $p$ in variables $X$ in a theory
$T$ is realized in the factor-algebra $\T _{\L}(X)/\theta_p$ at
the generic point  $\bar{x}=(x_1/\theta_p,\ldots,x_n/\theta_p)$,
so
 $$
\atp^{\T _{\L}(X)/\theta_p}(\bar{x}) \: = \: p.
$$
Indeed, for any atomic formula $t_1=t_2$, where $t_1,t_2 \in \Tr
_{\L}(X)$ one has  $(t_1=t_2)\in p$ if and only if  $t_1
\sim_{\theta_p} t_2$, which is equivalent to the condition  $\T
_{\L}(X)/\theta_p \models (t_1=t_2)$ under the interpretation
$x_i\mapsto x_i/\theta_p$.  The  generic point
$(x_1/\theta_p,\ldots,x_n/\theta_p)$ satisfies the following
universal property. If $p$ is realized in some $\L$-structure $\M$
at $(m_1, \ldots,m_n) \in M^n$ then the map $x_1 \rightarrow m_1,
\ldots, x_n \to m_n$ extends to a homomorphism  $\T
_{\L}(X)/\theta_p  \to \M$.

\begin{lemma}\label{type2}
Let $T$ be a universally axiomatized theory in $\L$. Then for any
finitely generated $\L$-structure $\M$ the following conditions
are equivalent:
\begin{enumerate}
\item [1)] \label{t1} $\M \in \Mod (T)$;
\item [2)] \label{t2}  $\M = \T _{\L}(X)/\theta_p$ for some complete atomic type $p$ in $T$.
\end{enumerate}
\end{lemma}

\begin{proof}
Let  $X = \{x_1, \ldots,x_n\}$ be a finite set and $\langle X \mid
S\rangle$ a presentation of an $\L$-structure $\M$, i.e., $\M
\cong T_{\L}(X)/\theta_S$.   If $p=\atp^{\M}(\bar{x})$,
$\bar{x}=(x_1,\ldots,x_n)$ then  $[S]=p^{+}$ and $\T
_{\L}(X)/\theta_p \cong \T_{\L}(X)/\theta_S \cong\M$. Therefore,
1)  implies  2).

To prove the converse, let  $p$ be an atomic type in $T$. We need
to show that $\T _{\L}(X)/\theta_p \in \Mod (T)$. Since $p$ is a
type in $T$ there exists a model  $\N \in \Mod (T)$ and a tuple of
elements $\bar{y} = (y_1,\ldots,y_n)\in N^n$ such that
$p=\atp^{\N}(\bar{y})$. If $\N'$ is a substructure of $\N$
generated by $y_1,\ldots,y_n$ then $\T _{\L}(X)/\theta_p \cong
\N'$. Since the theory $T$ is axiomatized by a set of universal
sentences  one has $\N' \in \Mod (T)$. Hence, $\T
_{\L}(X)/\theta_p \in \Mod (T)$.
\end{proof}


\subsection{Equationally Noetherian algebras}

The notion of equationally Noetherian groups was introduced in
\cite{BMR1}  and  \cite{BMRom}.

Let $\B$ be an algebra. For every natural number $n$  we consider
Zariski topology  on $B^n$.

 A subset  $Y\subseteq B^n$ is called {\em reducible} if it is
 a union of two proper closed subsets, otherwise, it is called  {\em irreducible}.

It is not hard to see that
an algebraic set $Y\subseteq B^n$ is irreducible if and only if it is not a finite union of proper algebraic subsets.

Recall, that a topological space is called {\em Noetherian} if it
satisfies the descending chain condition on closed subsets.

\begin{remark} \label{re:Noeth}
Let  $(W,\TT)$ be a topological space, $\AT$ a prebase of closed
subsets of  $\TT$, and $\BT$ the base of  $\T$, formed by the
finite unions of sets from $\AT$. Suppose that $\AT$ is closed
under finite intersections. Then the following conditions are
equivalent:
\begin{itemize}
\item the topological space  $(W,\TT)$ is Noetherian;
\item $\AT$ satisfies the descending chain condition.
\end{itemize}
In this case
\begin{enumerate}
\item [1)] the base  $\BT$ contains all closed sets in the topology  $\TT$;
\item [2)] any closed set  $Y$ in $\TT$ is a finite union of irreducible
closed sets from  $\AT$ ({\em irreducible components}): $Y=Y_1
\cup \ldots \cup Y_m$. Moreover, if $Y_i \not \subseteq Y_j$ for
$i \ne j$ then this decomposition is unique up to a permutation of
components.
\end{enumerate}
\end{remark}

\begin{definition} [No coefficients]
An algebra  $\B$ is \textit{equationally Noetherian}, if for any
natural number $n$ and any system of equations $S\subseteq
\At_{\L}(x_1,\ldots,x_n)$ there exists a finite subsystem $S_0
\subseteq S$ such that $\V _\B (S)=\V _\B (S_0)$.
\end{definition}

\begin{definition} [Coefficients in $\A$]
An $\A$-algebra  $\B$ is \textit{$\A$-equationally Noetherian} if
for any natural number $n$ and any system of equations $S\subseteq
\At_{\L_\A}(x_1,\ldots,x_n)$ there exists a finite subsystem $S_0
\subseteq S$ such that $\V _\B (S)=\V _\B (S_0)$.
\end{definition}

\begin{lemma}\label{noe}
An ($\A$-) algebra  $\B$ is ($\A$-) equationally Noetherian if and
only if for any natural number  $n$ Zariski topology on  $B^n$ is
Noetherian.
\end{lemma}

\begin{proof}
We prove the lemma for coefficient-free equations, a similar
argument  gives the result for equations with coefficients in
$\A$.

Assume  $\B$ is equationally Noetherian and consider a descending
chain of closed subsets  $Y_1 \supseteq Y_2 \supseteq Y_3
\supseteq \ldots$  of algebraic sets in  $\B^n$. Taking the
radicals one gets an ascending chain of subalgebras  $\Rad (Y_1)
\subseteq \Rad (Y_2) \subseteq \Rad (Y_3) \subseteq \ldots$. Put
$S=\bigcup_i {\Rad (Y_i)}$. By our assumption  the system $S$ is
equivalent to some finite subsystem $S_0 \subseteq S$. Clearly,
$S_0 \subseteq \Rad (Y_i)$ for some index $i$. Therefore, the
chains before stabilize.

Suppose now that for any natural number  $n$ Zariski topology on
$B^n$ is Noetherian. Let  $S\subseteq \At_\L (x_1,\ldots,x_n)$ be
an arbitrary system of equations in variables
$\{x_1,\ldots,x_n\}$.  Let  $(t_1=s_1) \in S$. If
$\V_\B(S)=\V_\B(\{t_1=s_1\})$ then there is nothing to prove.
Otherwise, there is an atomic formula $(t_2=s_2) \in S\backslash
\{t_1=s_1\}$ with $\V (\{t_1=s_1\})\varsupsetneq \V (\{t_1=s_1,
t_2=s_2\})$. Repeating this process  one can  produce a descending
chain of closed subsets in $B^n$.  Since $B^n$ is Noetherian the
chain is finite, so  $\V _\B (S)=\V _\B (S_0)$ for some finite
subsystem $S_0$ of $S$.
\end{proof}

The following result follows immediately from Lemma~\ref{noe} and
Remark~\ref{re:Noeth}.

\begin{theorem}\label{irr}
Let $\B$ be an ($\A$-)equationally Noetherian ($\A$-)algebra. Then
any algebraic set  $Y\subseteq B^n$ is a finite union of
irreducible algebraic sets ({\it irreducible components}): $Y=Y_1
\cup \ldots \cup Y_m$. Moreover, if  $Y_i \not \subseteq Y_j$ for
$i\ne j$ then this decomposition is unique up to a permutation of
components.
\end{theorem}

Now we give a characterization of the coordinate algebras of irreducible algebraic
sets over an arbitrary algebra $\B$.

\begin{lemma}\label{dis1}
Let $Y$ be an irreducible algebraic set over $\B$.  Then the
coordinate algebra  $\Gamma (Y)$ is discriminated by  $\B$.
\end{lemma}

\begin{proof}
Indeed, let $Y=\V (S)$ and  $\Gamma (Y)= \T_{\L}(X)/\theta_{\Rad
(Y)}$. Suppose, to the contrary, that there exist such atomic
formulas $(t_i=s_i)\in \At_{\L} (X)$, $(t_i=s_i)\not\in \Rad(Y)$,
$i=1, \ldots,m$, such that for any homomorphism   $h:\Gamma (Y)\to
\B$ there exists an index $i\in \{1,\ldots,m\}$ for which
$h(t_i/\theta_{\Rad (S)})=h(s_i/\theta_{\Rad (S)})$. This implies
that for any $p\in Y$ there exists an index  $i\in \{1,\ldots,m\}$
with $t_i^{\B}(p)=s_i^{\B}(p)$. Put $Y_i = \V (S\cup
\{t_i=s_i\})$, $i=1, \ldots,m$. Then  $Y=Y_1 \cup \ldots \cup Y_m$
and the sets $Y_1,\ldots,Y_m$ are proper closed subsets of $Y$~---
contradiction with irreducibility of $Y$. This shows that $\Gamma
(Y)$ is discriminated by  $\B$.
\end{proof}

The converse of this result also holds.

\begin{lemma}\label{dis2}
Let  $\Ce$ be a finitely generated  $\L$-algebra. If $\Ce$ is
discriminated by an $\L$-algebra  $\B$ then
 $\Ce$ is the coordinate algebra of some algebraic set over $\B$.
\end{lemma}

\begin{proof}
Since  $\Dis(\B)\subseteq \Res (\B)$ then by Proposition \ref{res}
$\Ce=\Gamma (Y)$ for some algebraic set $Y$ over $\B$. To prove
the result it suffices to reverse the argument in Lemma
\ref{dis1}. Indeed, suppose $Y=Y_1 \cup \ldots \cup Y_m$ for some
proper algebraic subsets $Y_i$. From  $Y_i\subset Y$ and $Y_i \neq
Y$ follows that  $\Rad (Y) \subset \Rad (Y_i)$ and $\Rad (Y) \neq
\Rad (Y_i)$, so there exists an atomic formula $(t_i=s_i)\in \Rad
(Y_i) \backslash \Rad (Y)$, $i= 1, \ldots,m$. This implies that
there is no any homomorphism $h:\Gamma (Y)\to \B$ with
$h(t_i/\theta_{\Rad (Y)}) \ne h(s_i/\theta_{\Rad (Y)})$ for all
$i=1, \ldots,m$,~--- contradiction with  $\Ce \in \Dis (B)$.
\end{proof}

\begin{theorem} \label{th:dis2}
Let $\B$ be an $\L$-algebra and   $\Ce$
a finitely generated  $\L$-algebra. Then $\Ce$  is the coordinate
algebra of some irreducible algebraic set over $\B$ if and only if
$\Ce$ is discriminated by $\B$.
\end{theorem}
 \begin{proof}
 Follows immediately from Lemmas \ref{dis1} and \ref{dis2}, and Remark \ref{re:Noeth}.
 \end{proof}

A similar argument gives the result for  $\A$-algebras.

\begin{theorem} \label{th:dis2-A}
Let $\B$ be an  $\A$-algebra and
$\Ce$ a finitely generated  $\A$-algebra. Then $\Ce$  is the
coordinate algebra of some irreducible algebraic set over $\B$ if
and only if $\Ce$ is $\A$-discriminated by $\B$.
\end{theorem}


\section{Limit algebras}

\subsection{Direct systems of formulas and limit algebras}
\label{subsec:limit}

In  this section we discuss limit $\L$-algebras. We need the
following notation. For a formula  $\varphi \in \Phi_\L(X)$ and a
map  $\gamma : X\to X'$  from $X$ into a set of variables
$X^\prime$  by  $\varphi(\gamma(X))$ we denote the formula
obtained from $\varphi$ by the substitution $x \to \gamma (x)$ for
every  $x\in X$.

\begin{definition}
A triple  $\Lambda=(I, \varphi_i,\gamma_{ij})$ is called a
\textit{direct system of formulas} in  $\L$ if
\begin{enumerate}
\item $\langle I, \leqslant \rangle$ is a directed ordering;
\item for each  $i\in I$ there is a finite reduct $\L
_i$ of $\L$ and a finite set of variables $X_i$ such that
$\varphi_i$ is a consistent diagram-formula in $\L_i$ in variables
$X_i$;
\item $\gamma_{ij}: X_i \to X_j$ is a map defined for every pair of indices $i,j\in I$, $i\leqslant j$,  such that:
\begin{itemize}
\item $\gamma_{ii}$ is the identical map for every  $i\in
I$;
\item $\gamma _{jk}\circ \gamma_{ij}=\gamma_{ik}$ for every  $i,j,k \in I$, $i\leqslant j\leqslant k$;
\item all conjuncts of $\varphi_i(\gamma_{ij}(X_i))$ are also conjuncts of $\varphi _j
(X_j)$;
\end{itemize}
\item for any  $c\in \L$ there exists  $i\in I$ such that  $\varphi _i$
contains a   conjunct of the type   $x_i=c$, where $x_i\in X_i$;
\item for any functional symbol  $F\in \L$, any $i\in I$, and any
tuple of variables $(x_1,\ldots,x_{n_{F}}) \in X_i^{n_F}$ there is
$j\in I$, $i\leqslant j$ such that $\varphi _j$ contains a
conjunct  of the type
$F(\gamma_{ij}(x_1),\ldots,\gamma_{ij}(x_{n_{F}}))=x_j$, where
$x_j\in X_j$.
\end{enumerate}
\end{definition}

Let  $\Lambda=(I, \varphi_i,\gamma_{ij})$ be a direct system of
formulas in $\L$. Define a factor-set $L(\Lambda)=\{(x_i,i), \;
x_i\in X_i, \: i\in I\}/\equiv$, where
$$
(x_i,i)\equiv (x_j,j) \quad \Leftrightarrow \quad \exists \: k\in
I, \; i,j\leqslant k,\quad \gamma_{ik} (x_i)=\gamma_{jk}(x_j).
$$
By $\langle x, i\rangle$ we denote the equivalence class of an
element
 $(x,i)$, $x\in X_i$, $i\in I$, relative to $\equiv$.

We turn the set  $L(\Lambda)$ into an $\L$-algebra interpreting
constants and operations from $\L$ on $L(\Lambda)$ as follows:
\begin{enumerate}
\item if $c\in \L$ is a constant symbol then
$c^{L(\Lambda)}=\langle x_i,i \rangle$, where  $i\in I$ is an
arbitrary index such that the conjunction $\varphi _i$ contains an
atomic formula of the type $x_i=c$, with $x_i \in X_i$;
\item if $F\in \L$ is a symbol of operation and $\langle x_1,i_1 \rangle,\ldots,\langle x_{n_{F}},i_{n_{F}}
\rangle \in L(\Lambda)$ then  $F^{L(\Lambda)}(\langle x_1,i_1
\rangle,\ldots,\langle x_{n_{F}},i_{n_{F}} \rangle)=\langle x_j,j
\rangle$, where  $j\in I$,  $i_1,\ldots, i_{n_{F}}\leqslant j$ and
such that the conjunction $\varphi _j$ contains a conjunct
$F(\gamma_{i_1 j}(x_1),\ldots,\gamma_{i_{n_{F}}
j}(x_{n_{F}}))=x_j$, $x_j\in X_j$.
\end{enumerate}

\begin{lemma}
The constants $c^{L(\Lambda)}$ and operations $F^{L(\Lambda)}$ are
well-defined.
\end{lemma}

\begin{proof}
Let $c\in \L$ be a constant symbol from $\L$. Then, from the
definition of the direct system of formulas, there exists $i\in I$
such that $\varphi _i$ contains a conjunct $x_i=c$ for some
$x_i\in X_i$. Suppose that there exists another index $j\in I$ for
which $\varphi _j$ contains a conjunct $x_j=c$, $x_j\in X_j$.
Since $\leqslant$ is a direct order on $I$ then there exists $k\in
I$ such that $i,j \leqslant k$. The formula  $\varphi_k$ contains
conjuncts $\gamma _{ik}(x_i)=c$ and $\gamma _{jk}(x_j)=c$. Since
$\varphi_k$ is realizable one has  $\gamma _{ik}(x_i)=\gamma
_{jk}(x_j)$, so $(x_i,i)\equiv (x_j,j)$. This shows that
$c^{L(\Lambda)}$ is well-defined.

Let $F\in \L$ be a functional symbol and  $\langle x_1,i_1
\rangle,\ldots, \langle x_{n_{F}},i_{n_{F}} \rangle \in
L(\Lambda)$. There exists  $j_0\in I$ such that $i_1,\ldots,
i_{n_{F}}\leqslant j_0$, in particular, $\gamma_{i_1
j_0}(x_1),\ldots,\gamma_{i_{n_{F}} j_0}(x_{n_{F}}) \in X_{j_0}$.
Then by the definition of the direct system there exists  $j\in I$
such that  $j_0\leqslant j$ and $\varphi _j$ contains a conjunct
$F(\gamma_{j_0 j}(\gamma_{i_1 j_0}(x_1)),\ldots,\gamma_{j_0
j}(\gamma_{i_{n_{F}} j_0}(x_{n_{F}}))=x_j$, $x_j\in X_j$. Since
$\gamma_{j_0 j}(\gamma_{i_k j_0}(x_k))=\gamma_{i_k j}(x_k)$, $k=
1, \ldots, n_{F}$, and  $i_1,\ldots, i_{n_{F}}\leqslant j$ then
$F^{L(\Lambda)}$ is defined on $\langle x_1,i_1 \rangle,\ldots,
\langle x_{n_{F}},i_{n_{F}} \rangle$.

Suppose there exists another  $i\in I$ such that $i_1,\ldots,
i_{n_{F}}\leqslant i$ and $\varphi _i$ contains a  conjunct
$F(\gamma_{i_1 i}(x_1),\ldots,\gamma_{i_{n_{F}}
i}(x_{n_{F}}))=x_i$ for some  $x_i\in X_i$. Then there exists  $k
\in I$ such that  $i,j \leqslant k$  and  $\varphi _k$ contains
the conjuncts $F(\gamma_{i_1 k}(x_1),\ldots,\gamma_{i_{n_{F}}
k}(x_{n_{F}}))=\gamma _{jk}(x_j)$ and  $F(\gamma_{i_1
k}(x_1),\ldots,\gamma_{i_{n_{F}} k}(x_{n_{F}}))=\gamma
_{ik}(x_i)$.  The diagram-formula $\varphi _k$ is consistent,
hence $\gamma _{ik}(x_i)=\gamma _{jk}(x_j)$ (otherwise $\varphi
_k$ should contain $\gamma _{ik}(x_i)\neq \gamma _{jk}(x_j)$ which
is impossible), therefore $(x_i,i)\equiv (x_j,j)$.

It is left to show that the value $F^{L(\Lambda)}(\langle x_1,i_1
\rangle,\ldots,\langle x_{n_{F}},i_{n_{F}} \rangle)$ does not
depend on the representatives  $(x_k, i_k)$ in the equivalence
classes $\langle x_k, i_k \rangle$, $k=1, \ldots,n_{F}$. The
argument is similar to the one above and we omit it.
\end{proof}

\begin{definition}
Let  $\Lambda=(I, \varphi_i,\gamma_{ij})$ be a direct system of
formulas in $\L$. Then the set $L(\Lambda)$ with the constants
$c^{L(\Lambda)}$   and operations $F^{L(\Lambda)}$ defined above
for $c , F \in \L$ is an $\L$-structure termed the  limit algebra
of $\Lambda$ or a limit algebra in $\L$.
\end{definition}

\begin{lemma}
Let  $\Lambda=(I, \varphi_i,\gamma_{ij})$ be a direct system of
formulas in $\L$. Then all formulas $\varphi_i$, $i\in I$,  hold
in the limit algebra $L(\Lambda)$ under the interpretation
$x\mapsto \langle x,i\rangle$, $x\in X_i$, $i\in I$.
\end{lemma}
\begin{proof}
 The result follows directly from the construction of the limit algebra $L(\Lambda)$.
\end{proof}

\begin{lemma}\label{lim1}
Let  $\Lambda=(I, \varphi_i,\gamma_{ij})$ be a direct system of
formulas in $\L$. Suppose $\{\B_i, i\in I\}$ is a family of
$\L$-algebras such that the formula $\varphi_i$ can be realized in
$\B_i$, $i\in I$. Then there is an ultrafilter  $D$ over $I$ such
that the algebra $L(\Lambda)$ embeds into the ultraproduct
$\prod_{i\in I}{\B_i}/D$.
\end{lemma}

\begin{proof}
By the conditions of the lemma for every $i\in I$ the formula
$\varphi _i$ holds in $\B_i$ under some interpretation $h_i: X_i
\to  B_i$. Define a map
$$f_0 : \{(x,i), x\in X_i,
i\in I\} \to \prod_{i\in I}{\B_i}$$ such that $f_0(x,i) = b \in
\prod_{i\in I}{\B_i}$, where $b(j)=h_j(\gamma_{ij}(x))$, if
$i\leqslant j$, and $b(j)= h_i(x)$, if $i > j$. To define an
ultrafilter $D$ over $I$ put $J_i = \{j\in I, i \leqslant j\}$,
$i\in I$. Since $\langle I, \leqslant \rangle$ is a direct
ordering any finite  intersection of sets from $D_0 = \{J_i \mid i
\in I\}$ contains a set from $D_0$. Therefore, there is an
ultrafilter $D$ on $I$ such that $D_0 \subseteq D$. Now, we define
a map  $f: L(\Lambda) \to \prod_{i\in I}{\B_i}/D$ by the rule:
$f(\langle x,i \rangle)=f_0(x,i)/D$, $\langle x,i \rangle \in
L(\Lambda)$.

It is not hard to verify that the map $f$ is well-defined and is
an injective $\L$-homomorphism.
\end{proof}

\begin{definition}
Let $\Lambda=(I, \varphi_i,\gamma_{ij})$ be a direct system of
formulas in $\L$ and $\B$ and $\L$-algebra. If every $\varphi_i$
can be realized in $\B$ then the limit algebra  $L (\Lambda)$ is
called a limit algebra over  $\B$. In this case we denote $L
(\Lambda)$ by $\B_{\Lambda}$.
\end{definition}

\begin{corollary}\label{lim2}
Let $\Lambda=(I, \varphi_i,\gamma_{ij})$ be a direct system of
formulas in $\L$,  $\B$ an  $\L$-algebra, and $\B_{\Lambda}$ a
limit algebra over $\B$. Then there exists an ultrafilter $D$ over
$I$ such that the  limit algebra $\B_{\Lambda}$ embeds into the
ultrapower  $\B^I/D$ of  $\B$.
\end{corollary}

The next result explain why limit algebras over $\B$ have this
name.
\begin{lemma}
Let $\Ce$ be a limit algebra over $\B$. Then $\Ce$ (viewed as a
structure in the relational language $\L^{rel}$)  is the limit of
a direct system of local submodels of $\B$.

\end{lemma}
\begin{proof}
Straightforward.
\end{proof}

Let $X=\{x_b, b\in B\}$ be a set of variables indexed by elements
from $B$. Now let $I$ be the set of all pairs  $(\L',X')$, where
$\L'$ is a finite reduct of the language $\L$ and $X'$ a finite
subset of $X$. Denote by $\B'$ $\L'$-reduct of $\B$ and by
$\varphi_{(\L',X')}$ the conjunction of all formulas $\phi$ such
that (i) $\phi$ or $\neg \phi$ is in the core diagram
$\Diag_0(\B')$, (ii) $V(\phi) \subseteq X'$, (iii) $\B \models
\phi$ under the interpretation $x_b\mapsto b$, $x_b \in X'$.
Clearly, $\varphi_{(\L',X')}$ is a diagram-formula. Conversely,
every diagram-formula realizable in $\B$ can be obtained in the
form of $\varphi_{(\L',X')}$. Define $(\L',X')\leqslant
(\L'',X'')$ if and only if when $\L'\subseteq \L''$ and
$X'\subseteq X''$.

It is easy to see that $(I,\leqslant)$ is a direct ordering.
Define the maps $\gamma _{(\L',X'), (\L'',X'')}$,
$(\L',X')\leqslant (\L'',X'')$, as the identical maps, i.e.,
$\gamma _{(\L',X'), (\L'',X'')}(x_b)=x_b$ for all  $x_b\in X'$.
Straightforward verification shows that
\[
\Lambda^{\B}\;=\;(\{(\L',X')\},\: \varphi_{(\L',X')}, \: \gamma
_{(\L',X'), (\L'',X'')})
\]
is a direct system of formulas in $\L$.

\begin{lemma}\label{lim3}
Let $\B$ be an $\L$-algebra and $\Lambda^{\B}$ the direct system
defined above. Then $L(\Lambda^{\B}) \cong \B$.
\end{lemma}

\begin{proof}
Notice that for every  $b_1,b_2\in B$ and  $i,j\in I$ the equality
$(x_{b_1},i) \equiv (x_{b_2}, j)$ holds in the limit algebra
$L(\Lambda^{\B})$ if and only if $b_1=b_2$. Therefore, the map
$f:\B \to L(\Lambda^{\B})$, defined by  $f(b)=(x_b,i)$ for any
$i\in I$, is a bijection. It is easy to check that $f$ is an
$\L$-homomorphism.
\end{proof}

\begin{lemma}\label{lim4}
Let $\B$ and  $\Ce$ be $\L$-algebras. If $\Th _{\exists} (\B)
\supseteq \Th _{\exists} (\Ce)$ then  $\Ce$
 is isomorphic to some limit algebra over $\B$.
\end{lemma}

\begin{proof}
By Lemma \ref{lim3} $L(\Lambda^{\Ce}) \cong \Ce$. The inclusion
$\Th _{\exists} (\B) \supseteq \Th _{\exists} (\Ce)$ shows that
all diagram-formulas in the direct system $\Lambda^{\Ce}$ are
realizable in $\B$. Hence $L(\Lambda^{\Ce})$ is a limit algebra
over $\B$.
\end{proof}


\subsection{Limit $\A$-algebras}

In this section we discuss limit algebras in the category of
$\A$-algebras.

\begin{definition}
Let $\A$ be an $\L$-algebra and $\L_\A$ the language $\L$ with
constants from $\A$. If $\B$ is an  $\A$-algebra and $\Lambda$ a
direct system of formulas in $\L_\A$ then the algebra
$\B_{\Lambda}$ is called a limit $\A$-algebra over $\B$.
\end{definition}

\begin{lemma}\label{lim5}
Let $\A$ be an $\L$-algebra, $\B$ an  $\A$-algebra,  and $\Lambda$
a direct system of formulas in $\L_A$. Then the limit algebra
$\B_{\Lambda}$ is an $\A$-algebra, i.e., $\B_{\Lambda} \models
\Diag (\A)$.
\end{lemma}

\begin{proof} It is not hard to prove the result directly from
definitions. However, it follows immediately from Corollary
\ref{lim2}.

\end{proof}

Since, in the notation above,  the limit algebra  $\B_{\Lambda}$
is an $\A$-algebra, all the results from Section
\ref{subsec:limit} hold (after an obvious adjustment) in the
category of $\A$-algebras. We just mention these results without
proofs.

\begin{corollary}
Let $\B$ be an $\A$-algebra in the language $\L_{\A}$ and
$\B_{\Lambda}$ the limit $\A$-algebra over  $\B$ relative to the
direct system $\Lambda=(I, \varphi_i,\gamma_{ij})$. Then there
exists an ultrafilter $D$ over  $I$ such that $\B_{\Lambda}$
$\A$-embeds into the ultrapower $\B^I/D$ of the algebra $\B$.
\end{corollary}

\begin{corollary}\label{lim6}
Let $\B$ be an $\A$-algebra in $\L_{\A}$ and $\Lambda_{\A}^{\B}$
be a direct system of formulas in $\L_A$, corresponding to $\B$
(see Lemma \ref{lim3}). Then  $L(\Lambda_{\A}^{\B}) \cong_{\A}
\B$.
\end{corollary}

\begin{corollary}\label{lim7}
Let $\A$ be an $\L$-algebra, $\B$ and  $\Ce$ $\A$-algebras  and
$\Th _{\exists, \A} (\B) \supseteq \Th _{\exists, \A} (\Ce)$. Then
$\Ce$ is $\A$-isomorphic to some $\A$-algebra which is a limit
algebra over $\B$.
\end{corollary}


\section{Unification Theorems}
\label{subsec:proof-MT} \label{sec6}

\medskip
\noindent {\bf Theorem A} {\it [No coefficients] Let $\B$ be an
equationally Noetherian algebra in a functional language $\L$.
Then for a finitely generated algebra  $\Ce$ of $\L$ the following
conditions are equivalent:
\begin{enumerate}

\item $\Th_{\forall} (\B) \subseteq \Th_{\forall} (\Ce)$, i.e., $\Ce \in \ucl(\B)$;
\item $\Th_{\exists} (\B) \supseteq
    \Th_{\exists} (\Ce)$;
\item $\Ce$ embeds into an ultrapower of $\B$;
\item $\Ce$ is discriminated by $\B$;
\item $\Ce$ is a limit algebra over $\B$;
\item $\Ce$ is  defined by a complete atomic type in the theory $\Th _{\forall} (\B)$ in
$\L$;
\item $\Ce$ is the coordinate algebra of an irreducible
algebraic set over $\B$ defined by a system of coefficient-free
equations.
\end{enumerate}
}

\medskip
\noindent {\bf Theorem B} {\it [With coefficients] Let $\A$ be an
algebra in a functional language $\L$ and $\B$ an
$\A$-equationally Noetherian $\A$-algebra. Then for a finitely
generated $\A$-algebra $\Ce$ the following conditions are
equivalent:
\begin{enumerate}

\item $\Th_{\forall,\A} (\B)
    \subseteq \Th _{\forall,\A} (\Ce)$, i.e., $\Ce \in \ucl_\A(\B)$;
\item $\Th_{\exists,\A} (\B) \supseteq \Th _{\exists,\A} (\Ce)$;
\item $\Ce$ $\A$-embeds into an  ultrapower of $\B$;
\item $\Ce$ is $\A$-discriminated by $\B$;
\item $\Ce$ is a limit  algebra   over  $\B$;
\item $\Ce$ is an algebra defined by a complete atomic type in
    the theory $\Th_{\forall,\A} (\B)$ in the language $\L_{\A}$;
\item $\Ce$ is the  coordinate algebra of an irreducible algebraic
set over $\B$ defined by  a system of equations with
coefficients in $\A$.
\end{enumerate}
}

\begin{proof}
We prove here only Theorem A, the argument for Theorem B is
similar and we omit it. Equivalence $1) \Longleftrightarrow 2)$ is
the standard result in mathematical logic.

Equivalence $1) \Longleftrightarrow 3)$ has been proven in  Lemma
\ref{le:ucl} (in the form $\ucl (\B)=\S\Pu(\B)$).

Equivalence $1) \Longleftrightarrow 6)$ has been proven in  Lemma
\ref{type2}.

To see that 1) is equivalent to 5) observe first that by Corollary
\ref{lim2} one has $5) \Longrightarrow 3)$, hence $5)
\Longrightarrow 1)$. The converse implication $1) \Longrightarrow
5)$  follows from Lemma \ref{lim4}.

Implication $4) \Longrightarrow 1)$ follows from $\Dis
(\B)\subseteq \ucl (\B)$ (see Lemma \ref{le:ucl}).

Now we prove the converse implication $1) \Longrightarrow 4)$.
Suppose that  $\Ce\not\in \Dis (\B)$. It suffices to show that
$\Ce \not\in \ucl (\B)$. Let $X=\{x_1,\ldots,x_n\}$ be a finite
set of generators of $\Ce$ and  $\langle X \mid S\rangle$ a
presentation of $\Ce$ in the generators $X$, where $S \subseteq
\At_{\L}(X)$. The latter means that $\Ce \simeq
\T_{\L}(X)/\theta_S$.

Since  $\B$ does not discriminate  $\Ce$ there are atomic formulas
$(t_i=s_i)\in \At_{\L} (X)$, $(t_i=s_i)\not\in [S]$, $i=1, \ldots,
m$, such that for any homomorphism $h:\Ce\to \B$ there is an index
$i\in \{1,\ldots,m\}$ for which $h(t_i/\theta_S)=h(s_i/\theta_S)$.
This means that for any point $p\in \V_\B(S)$ there is an index
$i\in \{1,\ldots,m\}$, with $t^{\B}_i(p)=s^{\B}_i(p)$. Since $\B$
is equationally Noetherian there exists a finite subsystem
$S_0\subseteq S$ such that $\V _\B(S_0) = \V _\B(S)$. Therefore,
the following universal statement holds in $\B$
\[
\forall \;  y_1 \ldots \forall \;  y_n \mathop \bigwedge
\limits_{(t=s) \in S_0} t (\bar{y})=s (\bar{y}) \; \to \; \bigvee
\limits_{i=1}^{m} t_i (\bar{x})=s_i (\bar{y}).
\]
On the other hand the formula
$$
\bigwedge \limits_{(t=s) \in S_0} t (\bar{y})=s (\bar{y}) \; \to
\; \bigvee \limits_{i=1}^{m} t_i (\bar{x})=s_i (\bar{y})
$$
is false in $\Ce$ under the interpretation $y_i\mapsto x_i$, $i=1,
\ldots,n$, hence   $\Ce \not\in \ucl (\B)$.

Equivalence  $4) \Longleftrightarrow  7)$ follows from Theorem
\ref{th:dis2}.

\end{proof}


\begin{remark}
In the case when $\A=\B$ the first two items in Theorem  B can be formulated in a more precise form:
$\Ce \equiv _{\forall,\A} \A$, and  $\Ce \equiv _{\exists,\A} \A$, correspondingly.
\end{remark}

\end{document}